\documentclass[12pt]{article}

\usepackage{enumitem}
\setlist[itemize]{label=\textbullet}

\usepackage[english]{babel}

\usepackage[utf8]{inputenc}
\usepackage[LGR, T1]{fontenc}
\usepackage{lmodern}
\usepackage[a4paper, margin=2cm, head=18.0pt]{geometry}
\usepackage{graphicx}
\usepackage{setspace}
\usepackage{booktabs}
\usepackage{tabularx}
\usepackage{multirow}
\usepackage{array}
\usepackage{caption}
\usepackage{float}
\usepackage{microtype}

\usepackage{mathtools}
\usepackage{amssymb}
\usepackage{amsthm}

\usepackage{xcolor}
\usepackage{tikz}
\usetikzlibrary{arrows.meta}
\usetikzlibrary{calc} % REQUIRED by the KAYROS wordmark coordinates below (paper0 gets it transitively via forest)
\usetikzlibrary{positioning}
\usetikzlibrary{shapes.geometric}
\usepackage{pgfplots}
\usepgfplotslibrary{groupplots}
\pgfplotsset{compat=1.9} % TeX Live 2016 ships pgfplots 1.13; 1.9 compat is safe

\definecolor{p5blue}{HTML}{3498DB}
\definecolor{p5red}{HTML}{E74C3C}
\definecolor{p5green}{HTML}{27AE60}
\definecolor{p5purple}{HTML}{9B59B6}
\definecolor{p5orange}{HTML}{E67E22}
\definecolor{p5grey}{HTML}{52514E}

\pgfplotsset{
  kayros axis/.style={
    axis lines=left,
    grid=both,
    grid style={p5grey!16, line width=0.3pt},
    tick label style={font=\footnotesize},
    label style={font=\small},
    title style={font=\normalsize, yshift=2pt},
    legend style={font=\footnotesize, draw=p5grey!45, fill=white, fill opacity=0.9,
                  text opacity=1, rounded corners=1.5pt},
    legend cell align=left,
    tick align=outside,
    tick style={p5grey!70, line width=0.5pt},
    every axis plot/.append style={line width=1.2pt},
  },
}
\tikzset{
  fignote/.style={font=\footnotesize},          % in-figure annotations, never on top of a curve
  figarrow/.style={-{Latex[length=2.6mm]}, line width=0.9pt, draw=p5grey!75},
}

\usepackage{csquotes}

\PassOptionsToPackage{hyphens}{url}
\usepackage[
  backend=biber,
  style=alphabetic,
  sorting=nyt,
  defernumbers=true,
  language=american,
  maxcitenames=1,
  mincitenames=1,
  maxbibnames=7,
  minbibnames=7
]{biblatex}
\DeclareSourcemap{
  \maps[datatype=bibtex]{
    \map{
      \step[fieldset=file, null]
    }
  }
}
\DeclareBibliographyAlias{softwareversion}{misc}
\DeclareBibliographyAlias{dataset}{misc}

\newcommand{\keywords}[1]{\textbf{Keywords}: \itshape #1 \normalfont}

\newcommand{\BoldSlashedO}{%
  \tikz[baseline=(char.base)]{
    \node[inner sep=0pt] (char) {O};
    \begin{scope}
    \clip (char.south west) rectangle (char.north east);
    \draw[line width=.55pt]
         ($(char.north west)+(4pt,-3pt)$) -- ($(char.south east)+(-4pt,3pt)$);
    \end{scope}
}}

\newcommand{\SlashedO}{%
  \tikz[baseline=(char.base)]{
    \node[inner sep=0pt] (char) {O};
    \begin{scope}
    \clip (char.south west) rectangle (char.north east);
    \draw[line width=.55pt]
         ($(char.north west)+(2.1pt,-1pt)$) -- ($(char.south east)+(-2.1pt,1pt)$);
    \end{scope}
}}

\DeclareRobustCommand{\KAYROSinner}{\textbf{KAYR\BoldSlashedO S}}
\DeclareRobustCommand{\kayrosinner}{KAYR\SlashedO S}
\newcommand{\KAYROS}{\texorpdfstring{\KAYROSinner}{KAYROS}}
\newcommand{\kayros}{\texorpdfstring{\kayrosinner}{KAYROS}}
\DeclareRobustCommand{\poryos}{\texorpdfstring{\textmd{\textsc{Poryos2026}}}{Poryos2026}}
\DeclareRobustCommand{\blauth}{\texorpdfstring{\textmd{\textsc{Blauth2024}}}{Blauth2024}}

\newcommand{\textgreek}[1]{\begingroup\fontencoding{LGR}\selectfont#1\endgroup}

\title{\Large \KAYROS: An Anytime and Exact Open-Source Solver for Duration-Minimization Time-Dependent Vehicle Routing\\[1.2ex]
\large A Technical Report and a Case Study in Human--AI Engineering}
\author{
  \normalsize Florian Rascoussier\textsuperscript{1,2}\thanks{ORCID identifier: \href{https://orcid.org/0009-0005-3253-9814}{0009-0005-3253-9814}.}\quad
  Romain Billot\textsuperscript{1}\thanks{ORCID identifier: \href{https://orcid.org/0000-0002-2614-3822}{0000-0002-2614-3822}.}\quad
  Lina Fahed\textsuperscript{1}\thanks{ORCID identifier: \href{https://orcid.org/0000-0002-7207-0580}{0000-0002-7207-0580}.}\quad
  Christine Solnon\textsuperscript{2}\thanks{ORCID identifier: \href{https://orcid.org/0000-0002-0919-496X}{0000-0002-0919-496X}.}\\[0.5ex]
  \small \textsuperscript{1}IMT Atlantique, Lab-STICC, CNRS, UMR 6285 (équipe DECIDE), France\\
  \small \textsuperscript{2}INSA Lyon, Inria, CITI, UR3720, 69621 Villeurbanne, France (équipe EMERAUDE)
}
\date{July 22, 2026. Revised October 19, 2026}

\usepackage[
    colorlinks=true,
    linkcolor=blue,
    citecolor=blue,
    urlcolor=blue
]{hyperref}
\hypersetup{
  pdftitle={KAYROS: An Anytime and Exact Open-Source Solver for Duration-Minimization Time-Dependent Vehicle Routing. A Technical Report and a Case Study in Human-AI Engineering},
  pdfauthor={Florian Rascoussier, Romain Billot, Lina Fahed, Christine Solnon}
}
\usepackage[capitalise, noabbrev]{cleveref}

\begin{document}

\maketitle
\pagenumbering{arabic}

% ########################################################
\begin{abstract}
Time-dependent routing recognizes that the same journey can take a different time depending on when it begins. Under duration minimization, even the departure time of each vehicle becomes a decision. Exact methods for this setting exist in the literature, but researchers and practitioners have lacked a ready-to-use open solver that combines rich piecewise-linear travel times, early feasible solutions and optimality claims. \kayros{} fills this gap with two modes on one checker-consistent engine: an Iterated Local Search that streams improving solutions and a Branch-Price-and-Cut method that can issue computational optimality certificates under explicit arithmetic and search assumptions. It installs with one command and has no proprietary dependency. The public MAMUT-routing store currently contains 704 \kayros{} certificates under a four-solve publication protocol, which has also led to the retraction and repair of invalid earlier claims. The report presents two complementary benchmark contributions to MAMUT-routing. The first integrates \blauth{}, a benchmark from the literature whose travel times derive from measured Uber speeds, for which \kayros{} provides new best-known solutions on all 40 instances. The second proposes \poryos{}, a new benchmark of 1,080 paired static and time-dependent instances built from OpenStreetMap road networks and controlled synthetic traffic. Finally, the report describes the intensive human--AI collaboration behind this work and the verification practices that kept its outputs independently verifiable.

\vspace{2ex}

\end{abstract}

\clearpage
\setcounter{tocdepth}{2}
\tableofcontents
\clearpage

% ########################################################
\section{Introduction}\label{sec:intro}
% Budget: ~1.5 pages total for 1.1-1.3.

\subsection{Context and claim}\label{sec:intro-claim}

\kayros{} \autocite{rascoussierKAYROS2026} is an open-source solver for the duration-minimization Time-Dependent Vehicle Routing Problem with Time Windows (TDVRPTW) and its time-window-free variant (TDVRP)\footnote{A.k.a TDCVRP, TDCVRPTW: capacity constraints included.}. It combines an anytime heuristic, which returns a feasible solution quickly and keeps improving it, with an exact Branch-Price-and-Cut (BPC) method that can establish computational optimality under the assumptions stated in \cref{sec:cert-semantics}. Both use the same piecewise-linear engine and report values recomputed by the public MAMUT-routing reference checker. The compatibility reference for this revision is \kayros{} 1.6.1 with \texttt{mamut-routing-lib} 0.11.0, available from PyPI and GitHub\footnote{\url{https://pypi.org/project/kayros/} and \url{https://github.com/0nyr/kayros}. Archived releases should be used when citing or reproducing a campaign.}.

To the best of our knowledge, \kayros{} is the first readily available solver to combine exact and anytime duration minimization for TDVRP(TW), rich piecewise-linear travel times rather than a fixed time discretization, one-command installation and a fully open dependency stack. Its function-composition engine follows \textcite{visserEfficientMoveEvaluations2020}, its local-search structures build on \textcite{blauthVehicleRoutingTimedependent2024}, and its exact component extends the open solver of \textcite{lera-romeroLinearEdgeCosts2020}. The claim concerns availability and scope, not a cross-solver performance ranking. The comparison campaign needed to make such a ranking responsibly is reserved for the lead author's thesis and the full solver paper.

The report's second contribution is a benchmark backstory that the solver paper cannot develop in the same detail. Two families distributed through MAMUT-routing \autocite{pichon:hal-05629810v1} provide complementary views of realistic TD routing. \blauth{} preserves delivery instances and high-effort references derived from measured Uber speeds \autocite{blauthVehicleRoutingTimedependent2024}. \poryos{}, designed and generated by the lead author, combines OpenStreetMap road networks for five cities with controlled demands, time windows and congestion. The first family emphasizes measured traffic and a delivery objective; the second emphasizes paired experimental control and reproducibility. \Cref{sec:poryos} explains why both are needed.

\subsection{Why this report also documents a collaboration}\label{sec:intro-collab}

This report also documents the intensive human--AI collaboration through which \kayros{} was built. \Cref{sec:collab} focuses on the reusable engineering method: how context was curated, how decisions and failures were recorded, how autonomous implementation was reviewed, and why external evidence remained necessary. Readers interested only in the technical contribution can skip that section. \Cref{sec:problem,sec:solver,sec:certification,sec:poryos,sec:results} are self-contained. Throughout, \enquote{the lead author} designates Florian Rascoussier, whose engineering work with AI tools is described in \cref{sec:collab}. The authorial \enquote{we} refers to the human coauthors.

\subsection{Scope and companion work}\label{sec:intro-scope}

This technical report favors explanation over formal development. \Cref{sec:problem} introduces TD routing through a two-customer route and documents the availability gap. \Cref{sec:solver} presents the solver at the architectural level, and \cref{sec:certification} explains what its certificates do and do not establish. \Cref{sec:poryos} develops the two benchmark families before \cref{sec:results} presents the public evidence. Formal models, proofs, algorithmic pseudocode, component ablations, composition-engine benchmarks and cross-solver results belong to the full paper and the thesis.

% ########################################################
\section{Time-dependent routing in a nutshell}\label{sec:problem}
% Budget: ~2.5 pages. Audience: NOT OR/VRP experts. Toy route <o,1,2,d> as running example.

This section assumes no vehicle-routing background: it states the problem, documents the availability gap, and builds every concept the rest of the report needs on a two-customer toy route.

\subsection{The problem and the gap}\label{sec:problem-gap}

The Vehicle Routing Problem with Time Windows (VRPTW) is one of the workhorses of Operations Research (OR): a fleet of capacitated vehicles based at a depot must visit a set of customers, each customer specifying a time window during which service may start, and the goal is to find feasible routes of minimum total cost \autocite{solomonAlgorithmsVehicleRouting1987}. The problem is NP-hard, ubiquitous in logistics, and has accumulated decades of competitive algorithmic literature.

The classic (static) model assumes that driving from $i$ to $j$ always takes the same time. Anyone who has commuted through a city knows this is false: the same road takes much longer at rush hours. The Time-Dependent VRPTW (TDVRPTW) drops this assumption and makes the travel time of every arc a function of the departure time \autocite{malandrakiTimeDependentVehicle1992}. The standard way to keep such a model physically sensible is the First-In First-Out (FIFO) property: departing later can never make you arrive strictly earlier \autocite{ichouaVehicleDispatchingTimedependent2003,fontaineExactAnytimeHeuristic2024}. Recent surveys of this field are given by \textcite{gendreauTimedependentRoutingProblems2015} and \textcite{adamoReviewRecentAdvances2024}.

\kayros{} targets the \emph{duration-minimization} variants of this problem, with or without time windows (TDVRPTW, TDVRP). Under this objective the departure time of each route from the depot is itself a decision variable, and a route's cost is the total time its vehicle is away, including any time spent waiting for a time window to open \autocite{visserEfficientMoveEvaluations2020, lera-romeroLinearEdgeCosts2020}. This is a natural objective (drivers are paid for their working time, not for the clock time at which they return), and a computationally demanding one: as \cref{sec:problem-toy} illustrates, even evaluating the cost of a \emph{single fixed route} requires optimizing over its departure time.

The research context also includes the time-dependent traveling-salesman work developed with Christine Solnon, from Pénélope Aguiar-Melgarejo's constraint-programming thesis \autocite{aguiar-melgarejoConstraintProgrammingApproach2016} to Romain Fontaine's thesis and the exact-and-anytime approach of Fontaine, Dibangoye and Solnon \autocite{fontaineExactAnytimeHeuristic2024,fontaineExactAnytimeApproach2023}. This line of work provides the single-vehicle foundations for the multi-vehicle extension addressed here.

Exact algorithms for these variants exist in the academic literature, notably the branch-and-price line of \textcite{dabiaBranchPriceTimeDependent2013} and the branch-price-and-cut (BPC) solver of \textcite{lera-romeroLinearEdgeCosts2020}, alongside exact approaches on related TD problems \autocite{ariglianoTimedependentAsymmetricTraveling2019, vuDynamicDiscretizationDiscovery2020, fontaineExactAnytimeApproach2023}. What has been missing is availability: to the best of our knowledge, no previously released solver combines, in one open-source package that installs with one command, (i) exact solving of the duration-minimization TDVRP(TW) over rich piecewise-linear travel-time functions rather than a time discretization, (ii) anytime behavior, meaning that the solver streams improving solutions from the first seconds and can be interrupted at any point with an honest answer, and (iii) no proprietary dependency, in particular no commercial LP solver. That combination is the claim this report makes for \kayros{}. It does not imply that other solvers cannot be adapted to related TD models, nor does this report use the later cross-solver campaign to claim a performance ranking.

\subsection{Why this gap matters: a decade of unmet demand}\label{sec:problem-demand}

The gap is practical, not merely bibliographic. For more than a decade, users have asked mainstream routing engines how to make travel times depend on the time of day. Their public issue trackers record the recurring answer: static matrices are central to fast move evaluation, existing extension points are partial, and a complete TD example or implementation is still missing. OR-Tools has received such requests since 2017, jsprit since 2013, VROOM explicitly relates its static matrix to constant-time move evaluation, and PyVRP tracks piecewise-linear arrival-time matrices as a future feature awaiting support\footnote{Representative threads, all accessed 2026-07-19: \href{https://github.com/google/or-tools/issues/339}{OR-Tools issue 339} and \href{https://github.com/google/or-tools/discussions/4230}{discussion 4230}, \href{https://github.com/PyVRP/PyVRP/issues/867}{PyVRP issue 867}, \href{https://github.com/VROOM-Project/vroom/issues/1280}{VROOM issue 1280}, and \href{https://github.com/graphhopper/jsprit/issues/22}{jsprit issue 22}.}. This sustained demand matters because adding a time index to a static matrix does not by itself solve the departure-time optimization illustrated below.

The closest available systems each provide part of the answer. OR-Tools exposes an unfinished cumul-dependent transit interface, but its maintainers state that the implementation does not work and is to be removed\footnote{\url{https://github.com/google/or-tools/issues/5070}.}; jsprit and Timefold allow custom time-dependent evaluators but provide no ready-made duration-minimization TDVRP(TW) model \autocite{graphhopperJsprit2026,timefoldSolver2026}; and Hexaly can express rich external functions, while its public TD template uses a small set of constant matrices and a different objective \autocite{hexalyTDCVRPTW2026}. These are useful capabilities, not interchangeable implementations of the problem studied here. A fair comparison consequently requires solver-specific models, disclosed conversions and independent route re-evaluation. That work was completed later and is reserved for the thesis and full paper. The narrower contribution of this report is the usable open combination stated in \cref{sec:intro-claim}.

\subsection{Travel times as functions: a two-customer route}\label{sec:problem-toy}

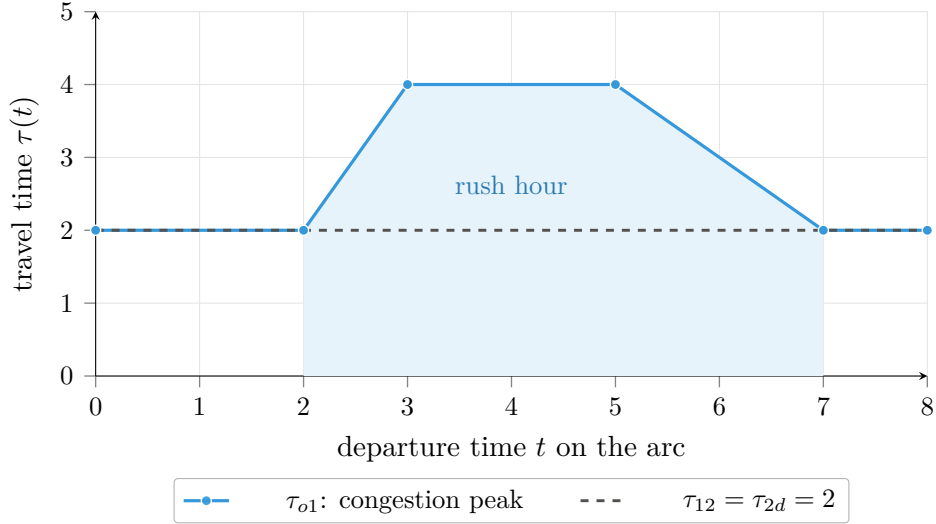
\begin{figure}[H]
\centering
% Toy example: the three arc Travel Time Functions (TTFs) of the route <o,1,2,d>.
% Exact coordinates verified against mamut-routing-lib (see scripts/toy_verify.py).
% Style: house palette and type sizes defined in main.tex (kayros axis).
\begin{tikzpicture}
\begin{axis}[
  kayros axis,
  width=0.74\textwidth, height=6.4cm,
  xlabel={departure time $t$ on the arc},
  ylabel={travel time $\tau(t)$},
  xmin=0, xmax=8, ymin=0, ymax=5,
  xtick={0,...,8}, ytick={0,...,5},
  legend style={at={(0.5,-0.28)}, anchor=north, legend columns=2, column sep=1.6em},
]
% shade the peak (drawn first so both line plots stay visible on top of the fill)
\addplot[fill=p5blue!12, draw=none, forget plot] coordinates {(2,2) (3,4) (5,4) (7,2)} \closedcycle;
% congestion peak on o->1
\addplot[p5blue, mark=*, mark size=1.9pt,
         mark options={fill=p5blue, draw=white, line width=0.5pt}]
  coordinates {(0,2) (2,2) (3,4) (5,4) (7,2) (8,2)};
\addlegendentry{$\tau_{o1}$: congestion peak}
% flat arcs
\addplot[p5grey, dashed, mark=none, line width=1.1pt] coordinates {(0,2) (8,2)};
\addlegendentry{$\tau_{12} = \tau_{2d} = 2$}
\node[fignote, p5blue!80!black] at (axis cs:4,2.62) {rush hour};
\end{axis}
\end{tikzpicture}
\caption{The three arc Travel Time Functions (TTFs) of the toy route $\langle o, 1, 2, d \rangle$. The arc from the depot to customer 1 has a morning congestion peak. The two other arcs are constant. Every travel time in this example is a continuous piecewise-linear function of the departure time on the arc.}
\label{fig:toy-ttf}
\end{figure}

Everything specific about time-dependent routing can be seen on a route with two customers, and this section walks through one. A vehicle leaves the depot $o$ at a departure time $t_0$ of our choosing, visits customer 1, then customer 2, and returns to the depot (written $d$ for the return copy of $o$, a standard convention). Customer 1 accepts service between times 4 and 9, customer 2 between 6 and 9, and service itself takes zero time in this toy. Time units are abstract. Think of hours of a working day. The travel time of each of the three arcs is a function $\tau(t)$ of the departure time $t$ on that arc: here the two arcs into and out of customer 2 always take 2 hours, while the arc from the depot to customer 1 takes 2 hours off-peak but climbs to 4 hours during a morning rush, as plotted in \cref{fig:toy-ttf}. These are Travel Time Functions (TTFs): piecewise-linear, and FIFO in the sense above.

Suppose the vehicle leaves at $t_0 = 0$. It arrives at customer 1 at time 2, but customer 1 only opens at 4: the vehicle waits two hours. It then reaches customer 2 at 6, right at opening, and returns to the depot at 8. The route takes 8 hours, 2 of them spent waiting. Now suppose it leaves at $t_0 = 2$ instead: it arrives at customer 1 exactly at 4, at customer 2 exactly at 6, and is back at 8, for a duration of 6 hours. Leaving two hours later gets the driver home at the same instant. Any later, though, and the departure runs into the rush hour on the first arc: arrival times now grow three times faster than the departure delay (the arc's travel time grows at slope 2, so arrival grows at slope 3), and past $t_0 = 3$ the vehicle can no longer reach customer 2 before its time window closes at 9: the route becomes infeasible.

\begin{figure}[H]
\centering
% Toy example: route ready-time function delta_r (left) and route duration function Delta_r (right)
% as functions of the depot departure time t0. Exact values verified with mamut-routing-lib:
% delta_r breakpoints (0,8),(2,8),(3,11); Delta_r(t0)=delta_r(t0)-t0; optimum 6 at t0=2; infeasible past 3.
% Style: house palette and type sizes defined in main.tex (kayros axis). Blue is the function
% being described, red the infeasible region, purple the EAT marker and green the optimum.
\begin{tikzpicture}
\begin{axis}[
  kayros axis,
  name=left,
  width=0.46\textwidth, height=6.4cm,
  xlabel={depot departure time $t_0$},
  ylabel={arrival back at depot $\delta_{\mathbf{r}}(t_0)$},
  xmin=0, xmax=4, ymin=5, ymax=12,
  xtick={0,1,2,3,4}, ytick={5,...,12},
]
% infeasible region (drawn first so the curve and its markers stay on top)
\addplot[fill=p5red!10, draw=none, forget plot] coordinates {(3,5) (4,5) (4,12) (3,12)} \closedcycle;
\addplot[p5blue, mark=*, mark size=1.9pt,
         mark options={fill=p5blue, draw=white, line width=0.5pt}]
  coordinates {(0,8) (2,8) (3,11)};
\node[fignote, p5red!75!black, rotate=90] at (axis cs:3.5,8.5) {infeasible};
\node[fignote, p5blue!80!black, anchor=north west] at (axis cs:0.1,7.86) {plateau: same arrival};
% EAT: Earliest Arrival Time = plateau value of delta_r
\addplot[mark=*, mark size=2.8pt, p5purple, only marks, forget plot] coordinates {(0,8)};
\node[fignote, p5purple, anchor=south west] at (axis cs:0.14,8.20) {EAT $= 8$};
\end{axis}
\begin{axis}[
  kayros axis,
  name=right,
  at={($(left.east)+(1.15cm,0)$)}, anchor=west,
  width=0.46\textwidth, height=6.4cm,
  xlabel={depot departure time $t_0$},
  ylabel={route duration $\Delta_{\mathbf{r}}(t_0)$},
  xmin=0, xmax=4, ymin=5, ymax=9,
  xtick={0,1,2,3,4}, ytick={5,...,9},
]
\addplot[fill=p5red!10, draw=none, forget plot] coordinates {(3,5) (4,5) (4,9) (3,9)} \closedcycle;
\addplot[p5blue, mark=*, mark size=1.9pt,
         mark options={fill=p5blue, draw=white, line width=0.5pt}]
  coordinates {(0,8) (2,6) (3,8)};
\node[fignote, p5red!75!black, rotate=90] at (axis cs:3.5,7.0) {infeasible};
% MDT: Minimum Duration Time = minimum duration value Delta* (y), at optimal departure t* (x)
\draw[densely dashed, p5green, line width=0.8pt] (axis cs:0,6) -- (axis cs:2,6);
\addplot[mark=*, mark size=3.2pt, p5green, only marks, forget plot] coordinates {(2,6)};
\node[fignote, p5green!75!black, anchor=south west] at (axis cs:0.14,6.10) {MDT $= 6$};
\node[fignote, p5green!75!black, anchor=north] at (axis cs:2.0,5.88) {$t^{*}_0 = 2$};
% The two slope labels ride on invisible paths laid along the segments they name, so `sloped`
% derives their angle from the canvas direction of that segment. Hand-set `rotate=` values do not
% survive: the angle on the page depends on the axis aspect ratio, so any change to the panel
% width, height or limits silently leaves the label off its line.
% Both are pulled off the midpoint, away from the vertex the two segments share, so the tail of
% one label and the head of the other do not crowd each other where the curve turns.
\draw[draw=none] (axis cs:0,8) -- (axis cs:2,6)
  node[fignote, p5grey, pos=0.43, above=2pt, sloped] {waiting-dominated};
\draw[draw=none] (axis cs:2,6) -- (axis cs:3,8)
  node[fignote, p5grey, pos=0.62, above=3pt, sloped] {congestion};
\end{axis}
\end{tikzpicture}
\caption{Ready time and duration of the toy route as functions of its departure time $t_0$. Departing at $t_0=2$ removes waiting and minimizes duration at $\Delta^{*}_{\mathbf{r}}=6$; departures after $t_0=3$ are infeasible. The reference checker computed all values.}
\label{fig:toy-duration}
\end{figure}
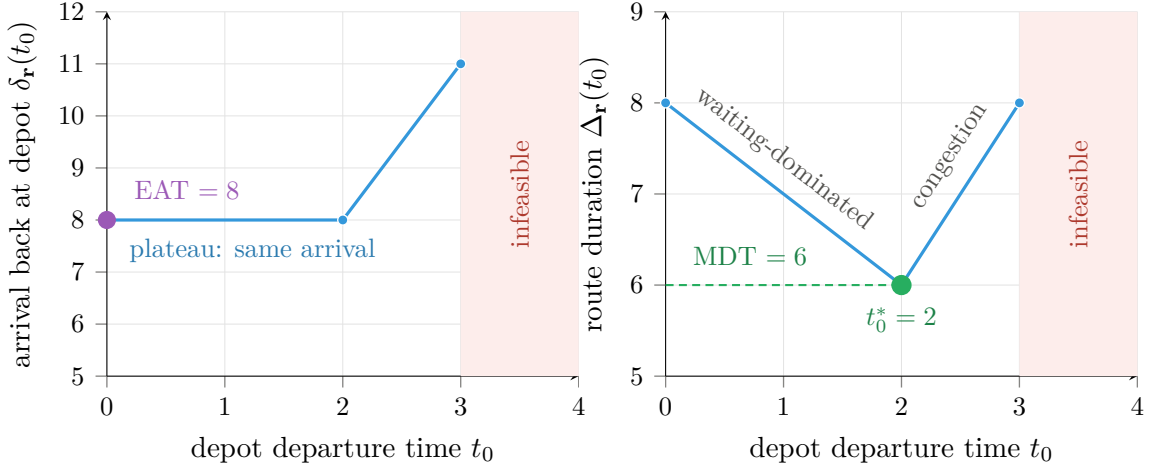

Two functions summarize this behavior completely, and they are the objects a time-dependent solver actually computes with. Composing the arc Arrival Time Functions (ATFs) $\alpha(t) = t + \tau(t)$ with the time-window waiting at each customer yields the route's \emph{ready-time function} $\delta_{\mathbf{r}}(t_0)$: the time the vehicle is back at the depot as a function of its departure time \autocite{visserEfficientMoveEvaluations2020, lera-romeroLinearEdgeCosts2020}. For our route, this entire analysis collapses to a three-breakpoint piecewise-linear function, plotted in \cref{fig:toy-duration} (left). Subtracting the departure time gives the route duration function $\Delta_{\mathbf{r}}(t_0) = \delta_{\mathbf{r}}(t_0) - t_0$, written $\Delta$ after \textcite{lera-romeroLinearEdgeCosts2020} (right panel): it falls at slope $-1$ while later departure eats into waiting, reaches the optimum $\Delta^{*}_{\mathbf{r}} = 6$ at $t_0 = 2$, then climbs at slope $+2$ into the congestion, until infeasibility at $t_0 = 3$. This minimum duration, $\Delta^{*}_{\mathbf{r}} = 6$, is the route's Minimum Duration Time (MDT). It is attained at the optimal departure time $t^{*}_0 = 2$, both marked on the right panel. Finding the best departure time of a route means minimizing this function.

The left panel also shows why the departure time is a genuine decision and not an afterthought. Every departure in $[0, 2]$ produces the \emph{same} arrival at 8, which is also the route's Earliest Arrival Time (EAT) at the depot: an objective that only minimizes the arrival time (the \emph{makespan}, natural when departures are fixed \autocite{fontaineExactAnytimeApproach2023}) cannot distinguish leaving at 0 from leaving at 2, while their durations differ by two full hours. Duration minimization sees the difference, which is precisely what makes it harder: for a route with many customers over travel-time functions with many pieces, $\delta_{\mathbf{r}}$ must be built by repeated function composition, and its number of breakpoints grows with both \autocite{visserEfficientMoveEvaluations2020}. A time-dependent solver performs this kind of composition millions of times: inside every local-search move evaluation of a heuristic and inside every labeling step of an exact pricing algorithm. Making that operation fast and \emph{binary64 checker-consistent}, meaning that it follows the reference checker's IEEE-754 binary64 evaluation rules and reproduces its results, is where much of the engineering of \cref{sec:solver} goes.

\subsection{Representing time dependence}\label{sec:problem-repr}

How the travel-time functions themselves are represented splits the field, and locating \kayros{} in that split is the last piece of context needed. The richest common representation, adopted by the exact literature \kayros{} builds on, keeps arc travel times as explicit continuous piecewise-linear functions, either given directly or generated from compact speed profiles in the model of \textcite{ichouaVehicleDispatchingTimedependent2003}, as in the canonical TD benchmark derived from Solomon instances by \textcite{dabiaBranchPriceTimeDependent2013}. \kayros{} consumes this representation natively: its engine composes Non-Decreasing Continuous PieceWise-Linear Functions (NDCPWLF, the ATFs of the previous section) using binary64 checker-consistent arithmetic.

The main alternative is to discretize time. Commercial solvers that handle time dependence typically model the horizon as a small number of constant-speed periods\footnote{For instance, the Hexaly modeling template for time-dependent routing divides the day into five parts (early morning, morning peak, day, evening peak, night): \url{https://www.hexaly.com/templates/time-dependent-routing-problem-with-time-windows-tdcvrptw}.}, and part of the exact literature works on time-expanded graphs whose discretization is refined on demand \autocite{vuDynamicDiscretizationDiscovery2020, heDynamicDiscretizationDiscovery2022}. Discretized models are legitimate and useful, but they answer a different computational question: a piecewise-constant approximation of a continuous travel-time landscape changes optimal values and can change optimal routes, so results on the two representations are not directly comparable. The claim of this report lives on the rich piecewise-linear side of this divide.

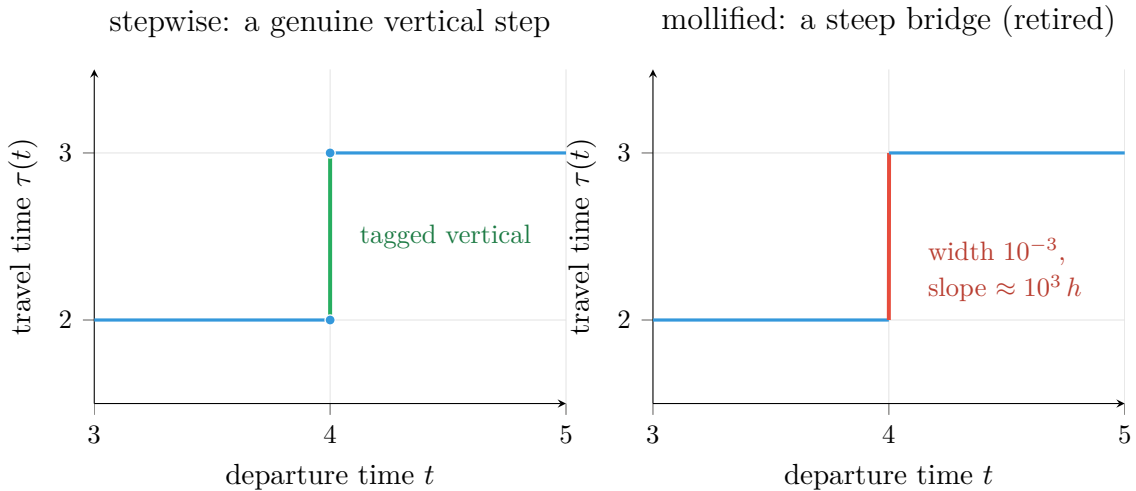
\begin{figure}[H]
\centering
% Stepwise travel times: a genuine value jump (left) versus the retired mollifier bridge (right).
% (caption text lives in main.tex)
% Style: house palette and type sizes defined in main.tex (kayros axis). Green marks the sound
% construct (the tagged vertical), red the retired unsound one (the steep bridge).
\begin{tikzpicture}
\begin{axis}[
  kayros axis,
  name=left,
  width=0.46\textwidth, height=6.0cm,
  xlabel={departure time $t$},
  ylabel={travel time $\tau(t)$},
  xmin=3, xmax=5, ymin=1.5, ymax=3.5,
  xtick={3,4,5}, ytick={2,3},
  title={stepwise: a genuine vertical step},
]
\addplot[p5blue, mark=none] coordinates {(3,2) (4,2)};
\addplot[p5blue, mark=none] coordinates {(4,3) (5,3)};
% the vertical, drawn as a first-class tagged object
\addplot[p5green, mark=none, line width=1.5pt] coordinates {(4,2) (4,3)};
\node[fignote, p5green!70!black, anchor=west] at (axis cs:4.08,2.5) {tagged vertical};
\addplot[mark=*, mark size=1.9pt, p5blue, only marks,
         mark options={fill=p5blue, draw=white, line width=0.5pt}] coordinates {(4,2) (4,3)};
\end{axis}
\begin{axis}[
  kayros axis,
  name=right,
  at={($(left.east)+(1.15cm,0)$)}, anchor=west,
  width=0.46\textwidth, height=6.0cm,
  xlabel={departure time $t$},
  ylabel={travel time $\tau(t)$},
  xmin=3, xmax=5, ymin=1.5, ymax=3.5,
  xtick={3,4,5}, ytick={2,3},
  title={mollified: a steep bridge (retired)},
]
\addplot[p5blue, mark=none] coordinates {(3,2) (4,2)};
\addplot[p5red, mark=none, line width=1.5pt] coordinates {(4,2) (4.001,3)};
\addplot[p5blue, mark=none] coordinates {(4.001,3) (5,3)};
\node[fignote, p5red!80!black, anchor=west, align=left]
  at (axis cs:4.12,2.30) {width $10^{-3}$,\\ slope $\approx 10^{3}\, h$};
\end{axis}
\end{tikzpicture}
\caption{A genuine travel-time jump (left) and the retired smoothing approach (right). The steep artificial bridge interacted unsafely with numerical tolerances; \kayros{} now represents the vertical step directly (\cref{sec:cert-refutation}).}
\label{fig:value-jump}
\end{figure}

One benchmark family deserves a special mention because it stress-tests that divide from within. The Lyon instances of \textcite{rifkiImpactSpatiotemporalGranularity2020}, derived from a real urban traffic simulation complemented with real-world data, have \emph{stepwise} travel times: at certain instants, the travel time of an arc genuinely jumps (\cref{fig:value-jump}, left). Such value jumps are vertical steps in the function graph, and they break the comfortable assumptions of continuous piecewise-linear machinery. How \kayros{} initially mishandled these jumps through a smoothing trick, how that error was caught by its own verification protocol, and how the engine now carries verticals as tagged first-class objects, is the story of \cref{sec:certification}. \Cref{sec:solver} first presents the solver in which that story takes place.

% ########################################################
\section{The \kayros{} solver}\label{sec:solver}
% Budget: ~3 pages + Figure 1 (architecture).

\kayros{} has two solving modes on one time-dependent engine, as sketched in \cref{fig:architecture}. The anytime mode produces a feasible solution early and improves it until the deadline. The exact mode searches for an optimum and reports whether it completed the conditions required for a computational certificate. Both send their final routes to an authority outside the solver: the MAMUT-routing reference checker \autocite{pichon:hal-05629810v1}. The name refers to the Greek concept of \emph{Kairos} (\textgreek{kair'os}), the opportune moment, because choosing \emph{when} a route departs is part of the problem. It is also a recursive acronym, \textbf{K}ayros \textbf{A}nytime-\textbf{Y}ielding \textbf{R}outing \textbf{O}ptimization \textbf{S}olver.

\begin{figure}[H]
\centering
\resizebox{0.97\textwidth}{!}{% KAYROS architecture: three layers inside the solver, the reference checker outside as referee.
% Style: house palette and type sizes defined in main.tex. Blue is the solver itself, red the
% external referee, grey the structural arrows and labels.
\begin{tikzpicture}[
  font=\small,
  layer/.style={draw=p5blue!65, line width=0.8pt, rounded corners=2.5pt, fill=p5blue!8,
                minimum width=10.4cm, minimum height=1.15cm, align=center},
  half/.style={draw=p5blue!65, line width=0.8pt, rounded corners=2.5pt, fill=p5blue!8,
               minimum width=5.05cm, minimum height=1.6cm, align=center},
  referee/.style={draw=p5red!75, line width=0.8pt, rounded corners=2.5pt, fill=p5red!8,
                  minimum width=2.7cm, minimum height=1.6cm, align=center},
  lbl/.style={font=\footnotesize\itshape, color=p5grey},
]
% API layer
\node[layer] (api) {Python API: \texttt{kayros.solve} (anytime), \texttt{kayros.lera} (exact) \quad \texttt{pip install kayros}};
% two middle boxes, centered under the API with a small fixed gap
% (anchoring them to api.south west/east let the wide API text push them apart
% and left no room for the arbitration label on the right)
\node[half, below=0.5cm of api.south, anchor=north east, xshift=-0.35cm] (anytime) {Anytime stack (C++)\\ greedy construction, TD-ILS\\ granular local search, TD-ACO};
\node[half, below=0.5cm of api.south, anchor=north west, xshift=0.35cm] (exact) {Exact component (C++)\\ vendored BPC\\ (Lera-Romero et al.\ 2020)\\ HiGHS LP backend};
% engine layer
\node[layer, below=3.15cm of api] (engine) {NDCPWLF engine (C++): checker-exact composition of travel-time functions,\\ exact IEEE-754 arithmetic, tagged verticals, bit-reproducible};
% referee outside
% the gap has to hold the "arbitration" label at readable type: at 1.9cm it overlapped
% both boxes once the figure was fitted to the text width.
\node[referee, right=2.8cm of exact.east] (checker) {MAMUT\\ reference\\ checker\\ (referee)};
% arrows
\draw[figarrow] (api.south -| anytime.north) -- (anytime.north);
\draw[figarrow] (api.south -| exact.north) -- (exact.north);
\draw[figarrow] (anytime.south) -- (anytime.south |- engine.north);
\draw[figarrow] (exact.south) -- (exact.south |- engine.north);
\draw[figarrow, <->, dashed, draw=p5red!80] (checker.west) -- node[above, font=\footnotesize, color=p5red!80!black] {arbitration} (exact.east);
\node[lbl, above=0.08cm of api.north west, anchor=south west] {solver};
\node[lbl, above=0.08cm of checker.north, anchor=south] {outside the solver};
\end{tikzpicture}}
\caption{The \kayros{} architecture. Two solving modes, anytime and exact, share one binary64 checker-consistent NDCPWLF engine. The checker version is pinned in \cref{sec:solver-packaging}. The checker itself stays outside the solver as the referee: reported route costs are recomputed by that checker.}
\label{fig:architecture}
\end{figure}

\subsection{Design principles}\label{sec:solver-principles}

Four commitments organize the implementation. First, the public checker is the referee: published costs come from an implementation outside \kayros{}, not from trusting the solver's internal objective. Second, both modes are genuinely anytime. They share one absolute deadline, provide improvements as they occur and return the latest checker-valid solution when interrupted. Third, termination is explicit: an exact run reports a computational optimum, an open search at the time limit, or a resource limit raised by the memory self-guard. Finally, the default installation is reproducible and fully open. Each run is single-threaded, the HiGHS Linear Programming (LP) solver is included in the wheels, and large campaigns obtain parallelism by launching independent runs rather than changing a run's arithmetic.

These commitments are safeguards, not proofs of correctness. In particular, reproducible arithmetic can reproduce the same defect on every machine. \Cref{sec:certification} explains why checker agreement, diversified solves and independent counterexamples must be used together.

\subsection{The NDCPWLF engine}\label{sec:solver-engine}

The bottom layer answers the question developed by the toy example: what does a fixed route cost when its departure time is optimized? It represents and composes Non-Decreasing Continuous PieceWise-Linear Functions (NDCPWLFs), following the move-evaluation results of \textcite{visserEfficientMoveEvaluations2020}. The C++ implementation follows the binary64 arithmetic of \texttt{mamut-routing-lib} 0.11.0, and an equivalence suite over the benchmark collection gates changes against that checker. Vertical pieces remain tagged because they can mean two different things: a genuine jump in travel time or a range of equally valid departure choices. Keeping that distinction is central to the repair discussed in \cref{sec:certification}.

\subsection{The anytime stack}\label{sec:solver-anytime}

The heuristic mode, exposed as \texttt{kayros.solve}, is a single-trajectory TD Iterated Local Search (ILS). A greedy constructor supplies the first solution; local search improves it through granular relocate, swap and 2-opt* moves; and ruin-and-recreate perturbations escape local optima under a late-acceptance rule \autocite{burkeLateAcceptanceHillClimbing2017}. Balanced route trees derived from \textcite{blauthVehicleRoutingTimedependent2024} make repeated TD move evaluations affordable. Every accepted move is then recomputed sequentially by the checker-consistent engine before it can become an incumbent. This separation lets fast data structures rank moves without allowing their evaluation order to redefine the published cost. A TD Ant Colony Optimization (ACO) strategy remains available, while the evidence in \cref{sec:results-ils-aco} motivates ILS as the default. The heuristic supports both total route duration and \texttt{FleetCostDuration}, which adds a fixed cost per vehicle.

\subsection{The exact component}\label{sec:solver-exact}

The exact mode, exposed as \texttt{kayros.lera}, extends the open BPC solver of \textcite{lera-romeroLinearEdgeCosts2020}. A master LP chooses routes from a growing pool, while a labeling algorithm searches for new routes that can improve that master solution. \kayros{} adds an open HiGHS backend, warm-start columns, deadline and memory-limit handling, TDVRP support, checker-consistent repricing and a dedicated path for stepwise travel times. These changes make the inherited research code usable as an interruptible component without hiding its assumptions.

The output is a computational certificate rather than a portable formal proof object. It establishes optimality under checker-consistent route costs, the stated LP and pricing tolerances, and the completeness of the labeling search modulo its epsilon-dominance rules. Safe, independently verified LP bounds would support a stronger claim, but they remain future work. For this reason, \kayros{} publishes certificates only after the redundant protocol of \cref{sec:cert-semantics} and retains explicit OPEN and RESOURCE\_LIMIT outcomes when that protocol does not finish.

\subsection{Packaging and availability}\label{sec:solver-packaging}

For this revision, install \kayros{} 1.6.1 with \texttt{mamut-routing-lib} 0.11.0 using \texttt{pip install "kayros==1.6.1" "mamut-routing-lib==0.11.0"} on supported Linux x86-64 Python versions. This pairing preserves the historical time-dependent checker contract. \kayros{} 1.6.1 excludes \texttt{mamut-routing-lib} 0.12, whose changed evaluation contract is incompatible with its internal evaluator. The earlier campaigns retain their original software provenance. The package includes both solving modes, pulls the benchmark loader and checker through \texttt{mamut-routing-lib}, and ships HiGHS without requiring a commercial license. The code is MIT-licensed, developed openly and archived by Software Heritage \autocite{rascoussierKAYROS2026}. Its public contract is deliberately small: an instance enters, and a stream of checker-priced solutions followed by an explicit verdict comes out.

% ########################################################
\section{Optimality certification, retraction and repair}\label{sec:certification}
% Budget: ~4 pages (was ~2 before the August 2026 incident subsection).

This section explains what \kayros{} means by a computational optimality certificate, how such claims are promoted, and why the public store records two retractions. The implementation details differ between the incidents, but the scientific lesson is common: agreement between several executions is useful evidence, yet it cannot establish correctness when every execution shares the same deterministic defect.

\subsection{What a certificate is}\label{sec:cert-semantics}

A \kayros{} certificate records that four exact solves reached the same checker-evaluated solution and completed the required exact-pricing phase. The four arms cross cold and warm starts with two labeling configurations. Every generated route entering the master problem is independently repriced, and a checker-infeasible route invalidates the run. The claim is therefore computational: the stored solution is optimal under the reference checker's route costs, the stated LP and pricing tolerances, and the completeness of the labeling search modulo epsilon dominance. It is not a proof object that can be checked without rerunning the solver, and it is not a claim of exact real arithmetic. It follows the practical tradition of the exact literature.

This definition also specifies honest failure. A run that reaches its deadline without completing the proof conditions returns OPEN; a run stopped by the memory self-guard returns RESOURCE\_LIMIT. Neither outcome is converted into a weaker form of certificate. The public metadata keeps the four run records and the software provenance needed to repeat the computation.

\subsection{The validation ladder}\label{sec:cert-ladder}

Before a build may publish certificates, it must pass four types of gate on Grid'5000. A cross-platform gate checks that independent builds reproduce the same route values. A differential fuzzer compares labeling modes on randomized instances. Targeted family sweeps exercise the arithmetic changed by the repair. Finally, the complete four-solve protocol is rerun over the affected certification domain. The criteria are fixed before execution: disagreement, checker-infeasible columns, incomplete exact pricing or an unexplained value change blocks promotion.

This ladder deliberately combines different failure detectors. Repetition tests stability, randomized comparisons expose configuration-dependent behavior, family sweeps cover representation-specific cases, and the external checker supplies counterexamples that do not depend on the prover. None is sufficient alone. After the second incident, the final gate covered all 1,352 instances of the four legacy families in both problem variants rather than only the family where the counterexample appeared.

\subsection{Retraction and repair of the stepwise certificates}\label{sec:cert-refutation}

The Rifki2020 family \autocite{rifkiImpactSpatiotemporalGranularity2020} contains genuine jumps in travel time. The inherited exact solver from Lera-Romero's PhD work was written for continuous functions, so the first extension replaced each vertical jump by a narrow, steep bridge. The transformation looked harmless pointwise, but it interacted with the pricing algorithm's numerical comparisons. Independent heuristic campaigns eventually produced checker-valid solutions better than 43 values that had been labelled optimal. All 160 certificates on the affected family were immediately retracted.

The repair replaced smoothing with a direct representation of value jumps. Vertical pieces now retain both endpoint values and a tag distinguishing a travel-time jump from a departure-time choice. The audit also found surrounding completeness defects in function inversion, route deduplication and pricing termination. The lasting protocol change is more important than that list: a certificate now requires an explicit audited exact-pricing phase, rather than inferring that exact pricing must have occurred from the solver's control flow.

This incident shows the useful side of diversification. Some invalid outcomes changed with the warm start or labeling configuration, while the independent heuristic supplied a decisive checker-valid counterexample. The repair was therefore tested from pinned counterexamples, differential fuzzing and a complete affected-family sweep before any certificate was restored. The numerical values and implementation chronology are retained in the public audit trail; the reusable lesson is that continuous approximations of stepwise semantics must not be treated as exact merely because they are narrow.

\subsection{A falsified certificate, and the re-derivation of the whole set}\label{sec:cert-recertification}

The second incident began when two heuristic configurations independently found a checker-valid solution below a published certificate on a continuous-travel-time Vu2020 instance. The exact build reproduced its earlier answer, so the disagreement was not random noise. At some breakpoints, a partial route can admit several valid durations because different departure choices lead to the same represented state. Pricing should retain the smallest of these durations, since it gives the most promising continuation. Instead, it evaluated only one representative value. When that value was too large, a genuinely improving route appeared unprofitable and was discarded. The same incomplete pricing produced the matching bound, so all four protocol arms agreed on the same wrong certificate.

The audit revealed a broader representation error. A vertical piece can describe either a discontinuous travel-time jump or the inverse of a flat interval containing several departure choices. Arithmetic designed to be safe for the first meaning is wrong for the second. The solver had selected between those semantics by inspecting intermediate operands, where the distinction was no longer recoverable. The repair now selects the arithmetic once from the instance's travel-time model and preserves the choice/jump distinction throughout the solve.

This case exposes the limit of redundant agreement more clearly than the stepwise incident. Four runs can detect a defect that depends on initialization or configuration, but four copies of one faulty reasoning path provide no protection against a shared deterministic error. The external checker verifies solution feasibility and cost, yet it cannot verify that the solver searched every improving column. A better checker-valid solution is therefore a decisive refutation, whereas the absence of such a solution is not a proof.

After the repair, the complete four-solve protocol was rerun from scratch over all 1,352 instances of the four legacy families in both problem variants. Of 468 pre-incident certificates, 466 were re-derived at their stored values. The falsified certificate and one further case with disagreement between protocol arms were withdrawn; their stored routes remain valid best-known solutions, but no optimality claim is attached to them. The same repaired build later produced the 238 \poryos{} certificates reported in \cref{sec:results-certificates}.

The durable conclusion is deliberately modest. A published \kayros{} certificate means that a diversified, checker-consistent computational protocol completed successfully under stated assumptions. It does not eliminate the possibility of a shared implementation defect, and the labeling arithmetic is not accompanied by independently verified safe dual bounds. This is why counterexample search continues after certification, why every contradiction triggers retraction, and why the public store preserves provenance rather than presenting certificate counts as immutable facts. This report preserves such provenance, which is often missing from the literature, and illustrates the difficulty of establishing correctness in a complex solver relying on subtle floating-point reasoning.

% ########################################################
\section{Complementary benchmark families: measured traffic and controlled variation}\label{sec:poryos}

A solver for realistic TD routing needs benchmarks that expose both the complexity of observed cities and the effect of controlled modeling choices. No single family provides both perfectly. This section presents two complementary contributions distributed through MAMUT-routing \autocite{pichon:hal-05629810v1}: \blauth{} preserves delivery instances derived from measured traffic speeds, while \poryos{} builds paired static and TD variants from retained road-network artifacts and reproducible synthetic traffic. The first supplies empirical travel-time structure and strong published references; the second makes it possible to vary congestion, time windows and problem class without changing the underlying customers.

The distinction also clarifies provenance. \blauth{} is a port of the benchmark of \textcite{blauthVehicleRoutingTimedependent2024}. Florian Rascoussier designed and generated \poryos{} \autocite{rascoussierPoryos2026}. MAMUT-routing provides the common schemas, checker, catalog and publication surface, while OpenStreetMap contributors provide the geographic source data used by both families \autocite{openstreetmapContributorsOpenStreetMap2026}.

\subsection{Measured traffic and delivery costs: \blauth{}}\label{sec:blauth}

\blauth{} contains 40 TDVRPTW instances across ten cities and four sizes from $n=10$ to $n=2000$ \autocite{blauthVehicleRoutingTimedependent2024}. Their arrival-time functions derive from Uber Movement speed observations and paths on OpenStreetMap road networks. The objective, \texttt{FleetCostDuration}, combines the optimal duration of each route with a fixed vehicle cost, reproducing the paper's model of \$200 per vehicle and \$20 per working hour. This changes the search problem substantially: saving one route can dominate a modest duration increase, so the heuristic must explore solutions that temporarily worsen time-window penalties or route duration before reducing the fleet.

The port preserves the published arrival-time functions, source hashes and initially the thirty high-effort BonnTour BKS at $n=500$, 1000 and 2000. The public checker re-evaluates each imported route under one objective convention. Preservation matters because the upstream traffic service has closed and the selected path geometries were not retained, so the original geographic derivation cannot now be reconstructed completely. Accordingly, \blauth{} offers measured-speed functions and strong references, but not the controlled pairing or full regeneration trail provided by \poryos{}.

\subsection{The paired \poryos{} design across four routing problems}\label{sec:poryos-design}

\poryos{} contains 1,080 instances derived from 60 bases: five cities, six customer counts $n \in \{10,25,50,100,500,1000\}$, and two sampling methods. The \emph{POI} method uses OpenStreetMap points of interest; the \emph{hybrid} method supplements them with a controlled spatial sample where POI coverage is sparse. Each base fixes the customers, coordinates, demands and vehicle capacity, with a capacity policy that prevents instances from degenerating into a single route.

\begin{table}[H]
\centering
% AUTO-GENERATED by scripts/gen_tables.py from the canonical Poryos2026 collection; do not edit by hand.
\begin{tabularx}{\textwidth}{lXrr}
\toprule
Problem type & Variants per base & Per base & Instances \\
\midrule
CVRP & Euclidean, shortest-road and fastest-road costs & 3 & 180 \\
VRPTW & Shared, tight and spread time-window sets & 3 & 180 \\
TDVRP & Two traffic models at three intensities & 6 & 360 \\
TDVRPTW & The same six overlays with shared audited windows & 6 & 360 \\
\midrule
\multicolumn{3}{l}{Total over 60 bases} & 1080 \\
\bottomrule
\end{tabularx}

\caption{Composition of \poryos{}. Every base shares its customer set, demands and capacity across its variants, and every instance carries a checker-valid best-known solution.}
\label{tab:poryos-design}
\end{table}

The three CVRP variants use Euclidean distance, shortest-road distance and free-flow fastest-road travel time. The three VRPTW variants use the fastest-road metric with shared, tight and spread time-window sets. The TDVRP and TDVRPTW sides each cross two traffic models with light, moderate and heavy intensities. The bare-base VRPTW instance and its six TDVRPTW twins share the same windows. The tight and spread sets are deliberately static-only. This construction supports controlled comparisons in which the customer geography stays fixed: Euclidean versus road-network cost, static versus time-dependent travel, windows versus no windows, and light versus heavy congestion.

\begin{figure}[H]
\centering
\resizebox{\textwidth}{!}{\input{figures/fig-poryos-map-n100}}
\caption{Road graph and checker-validated best-known solution for the TDVRPTW instance \texttt{poryos-lyon-n100-poi-bpr-heavy}. The 100 customers lie on the extracted Lyon road graph and the eight colored routes follow its directed edges. Geographic and road data are from OpenStreetMap contributors under ODbL 1.0 \autocite{openstreetmapContributorsOpenStreetMap2026}.}
\label{fig:poryos-map}
\end{figure}

Arc costs are stored as three-decimal values rather than rounded to the CVRPLIB integer convention, giving the static and time-dependent variants a common cost space. This choice is not cosmetic: scaling and rounding can alter feasibility, distort objectives and change empirical solver comparisons \autocite{rascoussier:hal-05646952}. Solvers that require integers can multiply all three-decimal quantities by 1,000 without discarding decimal precision, but published solutions are always re-evaluated by the MAMUT-routing checker in the canonical instance space. The checker, rather than a solver's internal arithmetic, determines feasibility and the stored BKS cost.

\subsection{From urban roads to reproducible time-dependent instances}\label{sec:poryos-generation}

\begin{figure}[H]
\centering
% Poryos2026 generation and validation pipeline, in natural reading order: stages 1 to 3 left to
% right on the top row, then 4 to 6 left to right on the bottom row, joined by a carriage-return
% arrow that runs the full width of the figure in the lane between the rows.
% The single row of six that this replaces had to be squeezed to \textwidth, which drove its
% labels down to 8pt on the page; two rows of three carry the same content at readable type.
% Style: house palette and type sizes defined in main.tex.
\begin{tikzpicture}[
  font=\small,
  box/.style={rectangle, rounded corners=2.5pt, line width=0.8pt, align=center,
              text width=4.1cm, minimum height=2.55cm, inner sep=6pt},
  data/.style={box, fill=p5blue!8, draw=p5blue!65},
  synth/.style={box, fill=p5orange!10, draw=p5orange!80},
  gate/.style={box, fill=p5green!9, draw=p5green!70},
  output/.style={box, fill=p5purple!8, draw=p5purple!65},
]
\node[data] (osm) {\textbf{1. City data}\\[3pt]OpenStreetMap city extract};
\node[data, right=0.85cm of osm] (graph) {\textbf{2. Road foundation}\\[3pt]Drivable graph, real geometry, class-derived free-flow speeds};
\node[synth, right=0.85cm of graph] (base) {\textbf{3. Customer base}\\[3pt]POI or hybrid sampling; demands, capacity and service times};
\node[synth, below=1.6cm of osm] (derive) {\textbf{4. Derivations}\\[3pt]Three static metrics, three window sets, six traffic overlays, pinned paths and ATFs};
\node[gate, right=0.85cm of derive] (audit) {\textbf{5. Validation}\\[3pt]Capacity and horizon anchor audit under all six overlays, SHA-256 pins, reference checker};
\node[output, right=0.85cm of audit] (publish) {\textbf{6. Published family}\\[3pt]CVRP $\times 3$, VRPTW $\times 3$, TDVRP $\times 6$, TDVRPTW $\times 6$, complete BKS};
\draw[figarrow] (osm) -- (graph);
\draw[figarrow] (graph) -- (base);
\draw[figarrow] (derive) -- (audit);
\draw[figarrow] (audit) -- (publish);
% Carriage return from stage 3 to stage 4: down out of the last box of the top row, straight back
% along the lane halfway between the rows, then down into the first box of the bottom row. The
% lane coordinate is derived from the boxes themselves, so the arrow stays centred in the gap if
% the row separation is ever changed.
\coordinate (lane) at ($(base.south)!0.5!(derive.north)$);
\draw[figarrow, rounded corners=6pt]
  (base.south) -- (base.south |- lane) -- (derive.north |- lane) -- (derive.north);
\end{tikzpicture}
\caption{Generation, feasibility and publication pipeline for \poryos{}, read left to right and top to bottom. Blue boxes are OSM-derived foundations, orange boxes are controlled benchmark derivations, the green box gathers validation gates, and the purple box summarizes the published paired variants.}
\label{fig:poryos-pipeline}
\end{figure}
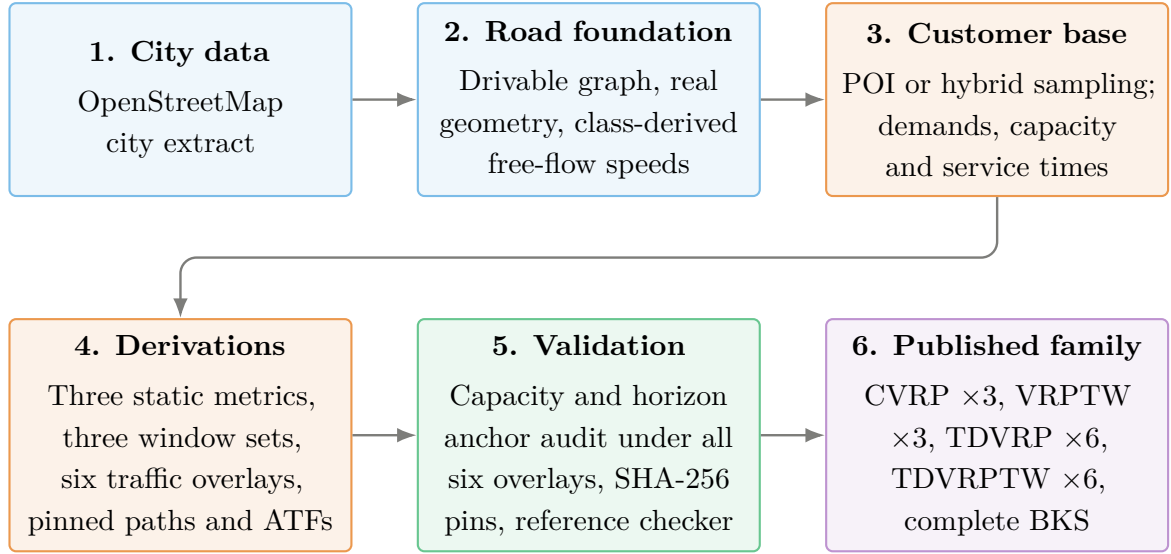

\poryos{} turns a city extract into a benchmark through a visible sequence of transformations, summarized in \cref{fig:poryos-pipeline}. The published collection includes the resulting instances and the artifacts needed to inspect them; the website publishes these canonical files rather than generating data on demand.

First, the generator extracts a directed drivable graph from OpenStreetMap, retains its geometry and road classes, samples customers, and computes the three static cost views. Each customer pair receives one pinned free-flow fastest path. Congestion changes the traversal time of that path but does not reroute the vehicle dynamically. This restriction makes static and TD variants comparable and leaves time-dependent path choice as a separate future dimension.

Second, the generator adds controlled routing data. Demands, capacities and service times are synthetic, and Solomon-like route-centered anchor solutions guide the time windows \textcite{solomonAlgorithmsVehicleRouting1987}. Two traffic models then create light, moderate and heavy congestion. One follows the Bureau of Public Roads volume-delay relationship from simulated commuter flows \autocite{federalHighwayAdministrationMultiresolutionModeling2022}; the other imposes reproducible morning and evening waves. Both slow the same free-flow network, so increasing congestion never creates an artificial speedup.

Finally, hourly edge speeds are composed along the pinned paths into FIFO arrival-time functions. Deterministic sampling and simplification keep the representation tractable, while hashes bind every instance to its graph, traffic overlay and materialized functions. A clone of the published collection is therefore sufficient to load and check the canonical benchmark without contacting a routing service or relying on the current state of OpenStreetMap.

\subsection{Feasibility by construction, then by audit}\label{sec:poryos-feasibility}

Random time windows and congestion can make large instances infeasible before a solver is ever tested. \poryos{} therefore carries a witness through generation. Deterministic anchor routes cover every customer, respect capacity and fit the horizon under free-flow travel. Before the TD twins are emitted, those routes are simulated under all six traffic scenarios. If congestion violates a deadline, the smallest shared deadline adjustment that restores the witness is applied consistently across the paired variants; any failure of coverage, capacity, positive window width or return before the horizon aborts generation.

The witness establishes existence, not quality. Independent validation then checks schemas and hashes, materializes the arrival-time functions and evaluates the solution files. Every published instance carries a checker-valid Best-Known Solution (BKS), while the 238 computational optima reported in \cref{sec:results-certificates} cover only the smaller TD instances. The checksummed published bytes define the benchmark. A fresh extraction from a later OpenStreetMap snapshot creates a new derivation rather than silently changing an existing instance.

% ########################################################
\section{Experimental evidence on MAMUT-routing}\label{sec:results}
% Budget: ~2.5 pages + Tables 1-2.

\subsection{The benchmark collection}\label{sec:results-collection}

\kayros{} is developed and evaluated against the open MAMUT-routing platform \autocite{pichon:hal-05629810v1}. The \href{https://mamut-routing.univ-ubs.fr/}{official website} and the lead author's \href{https://mamut-routing.onyr.net/}{mirror} publish the same catalog and checker-facing artifacts.\footnote{The mirror was established because the official host had been unavailable for more than a month at the time of writing.} The current catalog spans four routing problems, 18 families and 4,718 instances, each with a checker-valid BKS.

The TD collection joins several kinds of evidence. The legacy families of \textcite{dabiaBranchPriceTimeDependent2013}, \textcite{ariglianoTimedependentAsymmetricTraveling2019}, \textcite{vuDynamicDiscretizationDiscovery2020} and \textcite{rifkiImpactSpatiotemporalGranularity2020} support comparison with earlier research. Lera2026 and \poryos{} extend paired TDVRP and TDVRPTW coverage to $n=1000$, while \blauth{} reaches $n=2000$ under a delivery-specific objective. Their travel-time functions are published as checksummed artifacts, and their BKS were produced or audited through Grid'5000 campaigns. This shared checker space is what makes the following comparisons meaningful.

\subsection{Optimality certificates in the public store}\label{sec:results-certificates}

Under the protocol of \cref{sec:cert-semantics}, the public store contains \textbf{704 computational optimality certificates}\footnote{As of 2026-08-08}. The repaired build re-derived 466 certificates over the 1,352 legacy instances and produced 238 over the 720 TD \poryos{} instances. \Cref{tab:certificates} shows where they lie. Coverage reaches from $n=10$ to $n=100$, but it is uneven: 165 of 168 Vu2020 TDVRPTW instances certify, whereas removing time windows weakens pruning and drives more TDVRP runs to the memory frontier. On \poryos{}, almost all instances at $n\le25$ certify, while only six TDVRPTW instances do at $n=50$ and none above that size.

OPEN and RESOURCE\_LIMIT outcomes are part of this result. They mark the present reach of the exact component rather than being silently omitted failures. No certificate is claimed for Lera2026 or \blauth{}, whose sizes or objective lie outside that reach. Moreover, the 173 certificates that improved a previous BKS mostly improve solutions produced by \kayros{}' own heuristic. They measure the relationship between the two solver modes, not a literature comparison; the latter is restricted to the cases below.

\begin{table}[H]
\centering
% AUTO-GENERATED by scripts/gen_tables.py from td_recert3/verdicts_final.json (legacy families)
% and the public Poryos2026 optimality flags cross-checked against poryos_cert/verdicts_final.json;
% do not edit by hand.
\begin{tabular}{llrrr}
\toprule
Problem type & Family & Instances & Proven optimal & Share \\
\midrule
TDVRPTW & Dabia2013 & 168 & 114 & 68\% \\
TDVRPTW & Vu2020 & 168 & 165 & 98\% \\
TDVRPTW & Rifki2020 & 180 & 71 & 39\% \\
TDVRPTW & Ari2018 & 160 & 40 & 25\% \\
TDVRPTW & Poryos2026 & 360 & 123 & 34\% \\
\midrule
TDVRP & Dabia2013 & 168 & 40 & 24\% \\
TDVRP & Vu2020 & 168 & 0 & 0\% \\
TDVRP & Rifki2020 & 180 & 36 & 20\% \\
TDVRP & Ari2018 & 160 & 0 & 0\% \\
TDVRP & Poryos2026 & 360 & 115 & 32\% \\
\midrule
\multicolumn{2}{l}{Total (five TD families)} & 2072 & 704 & 34\% \\
\bottomrule
\end{tabular}

\caption{Computationally certified optima by family in the public MAMUT-routing store (2026-08-08), under the four-solve protocol of \cref{sec:cert-semantics}. Counts are generated from the campaign records and match the published per-instance metadata.}
\label{tab:certificates}
\end{table}

\subsection{Comparing against published results}\label{sec:results-literature}

The Dabia2013 TDVRPTW distribution is the strongest machine-checkable link with earlier exact work \autocite{lera-romeroLinearEdgeCosts2020}. Of 114 current \kayros{} certificates on that family, 103 reproduce a published optimum. Ten others establish values below references that the earlier tables left unproven, with improvements from 0.06\% to 2.15\%. One final case differs only through checker re-evaluation of the published route and remains below the reporting threshold. No \kayros{} certificate contradicts a published proven optimum. The ten substantive improvements are listed in \cref{tab:improvements}; the body retains the conclusion rather than the per-instance arithmetic.

The Vu2020 references are heuristic values rather than published optima. We therefore believe that the 165 TDVRPTW certificates and the 76 TDVRP certificates provide new exact results, but this has not been exhaustively reconciled with every original table of \textcite{vuDynamicDiscretizationDiscovery2020}. The qualification matters: a public certificate establishes the solver's stated computational claim, while priority over the literature is a separate bibliographic claim.

\blauth{} provides a complementary heuristic comparison under \texttt{FleetCostDuration}. Long-budget \kayros{} searches improve all thirty high-effort BonnTour references, but only narrowly: reductions range from 0.032\% to 1.314\%, with a mean of 0.254\%, and become smaller as size increases. Searches constrained to one fewer vehicle found no feasible reduction at the tested budgets. Hence the result supports two conclusions: \kayros{} is competitive on this family, and the original BonnTour solutions remain strong baselines. It is not a time-matched anytime solver comparison. \Cref{tab:blauth-improvements} preserves the per-instance values in the appendix.

\subsection{Complete heuristic coverage of \poryos}\label{sec:results-poryos}

The public store provides a checker-valid BKS for every \poryos{} instance, from $n=10$ to $n=1000$. \kayros{} produced 719 of the 720 TD records; one was subsequently improved by an external solver (Timefold). Among the 238 certified TD instances, the heuristic had already reached the certified value in 234 cases, while the exact component made four small corrections. Complete coverage makes the family usable immediately, but it must not be confused with complete optimization: above $n=50$, every TD record remains a best-found solution. Canonical checker evaluation is especially important here because scaling and rounding can change feasibility and empirical comparisons \autocite{rascoussier:hal-05646952}.

\subsection{Choosing the default heuristic}\label{sec:results-ils-aco}

The main internal algorithmic decision was choosing ILS over ACO as the default anytime strategy. A matched 20,808-run comparison covered five TD families and sizes from $n=10$ to $n=1000$. ILS won the large majority of paired cells, and its advantage increased with instance size. This is evidence about the two implementations on the present benchmark design, not a general theorem about ILS and ACO.

\begin{figure}[H]
\centering
% AUTO-GENERATED by scripts/gen_figures.py from convergence campaign CSVs; do not edit by hand.
\begin{tikzpicture}
\begin{groupplot}[kayros axis,group style={group size=3 by 2,horizontal sep=1.75cm,vertical sep=1.45cm},scale only axis,width=3.75cm,height=3.4cm,grid=major,minor tick style={draw=none},xmode=log,xmin=1,unbounded coords=jump]
\nextgroupplot[title={Dabia2013, $n=100$},ylabel={mean gap to BKS (\%)},ymin=0,xmax=1200]
\addplot[p5blue] coordinates {(1,6.7297) (2,5.2535) (3,4.5692) (5,3.8558) (7,3.4979) (10,3.1167) (15,2.7798) (20,2.5961) (30,2.3021) (45,2.0273) (60,1.8393) (90,1.584) (120,1.4071) (180,1.1319) (240,0.9736) (360,0.7743) (480,0.6867) (600,0.6489) (900,0.5913) (1200,0.5601)};
\addplot[p5orange,dash pattern=on 4.5pt off 2.5pt] coordinates {(1,7.7368) (2,7.2) (3,6.6832) (5,6.1764) (7,5.8173) (10,5.4936) (15,5.1382) (20,4.9954) (30,4.7729) (45,4.4974) (60,4.3712) (90,4.1856) (120,4.0407) (180,3.8031) (240,3.5251) (360,3.3444) (480,3.2895) (600,3.2607) (900,3.2483) (1200,3.2334)};
\nextgroupplot[title={within 1\%},xmin=1,xmax=1200,ymin=0,ymax=1,ytick={0,0.25,0.5,0.75,1},ylabel={fraction of runs}]
\addplot[p5blue] coordinates {(1,0.1929) (2,0.2143) (3,0.225) (5,0.2357) (7,0.25) (10,0.2643) (15,0.2786) (20,0.2821) (30,0.3071) (45,0.3357) (60,0.3607) (90,0.3857) (120,0.425) (180,0.5214) (240,0.5964) (360,0.7214) (480,0.7786) (600,0.7893) (900,0.8107) (1200,0.8393)};
\addplot[p5orange,dash pattern=on 4.5pt off 2.5pt] coordinates {(1,0.1607) (2,0.1607) (3,0.1607) (5,0.1643) (7,0.1679) (10,0.1786) (15,0.1893) (20,0.1929) (30,0.1929) (45,0.2036) (60,0.2071) (90,0.2107) (120,0.2143) (180,0.2179) (240,0.2179) (360,0.2321) (480,0.2321) (600,0.2357) (900,0.2429) (1200,0.2429)};
\nextgroupplot[title={within 0.5\%},xmin=1,xmax=1200,ymin=0,ymax=1,ytick={0,0.25,0.5,0.75,1}]
\addplot[p5blue] coordinates {(1,0.1714) (2,0.1857) (3,0.1964) (5,0.2107) (7,0.2143) (10,0.2393) (15,0.2393) (20,0.25) (30,0.25) (45,0.2607) (60,0.2714) (90,0.2821) (120,0.3107) (180,0.3643) (240,0.3786) (360,0.45) (480,0.5036) (600,0.5536) (900,0.5857) (1200,0.6071)};
\addplot[p5orange,dash pattern=on 4.5pt off 2.5pt] coordinates {(1,0.1607) (2,0.1607) (3,0.1607) (5,0.1607) (7,0.1643) (10,0.1679) (15,0.1679) (20,0.1679) (30,0.1714) (45,0.1821) (60,0.1821) (90,0.1857) (120,0.1893) (180,0.1929) (240,0.2143) (360,0.2179) (480,0.2214) (600,0.2214) (900,0.2214) (1200,0.2214)};
\nextgroupplot[title={Poryos2026, $n=1000$},xlabel={time (s)},ylabel={mean gap to BKS (\%)},ymin=0,xmax=7200]
\addplot[p5blue] coordinates {(90,2.0332) (120,1.8352) (180,1.6436) (240,1.4738) (360,1.1713) (480,0.9882) (600,0.8628) (900,0.6877) (1200,0.5827) (1800,0.449) (2400,0.3514) (3600,0.2286) (4800,0.1499) (7200,0.0579)};
\addplot[p5orange,dash pattern=on 4.5pt off 2.5pt] coordinates {(90,2.0079) (120,1.8373) (180,1.7372) (240,1.7372) (360,1.7372) (480,1.733) (600,1.7078) (900,1.6916) (1200,1.6742) (1800,1.6466) (2400,1.6237) (3600,1.588) (4800,1.5693) (7200,1.5382)};
\nextgroupplot[title={within 1\%},xlabel={time (s)},xmin=1,xmax=7200,ymin=0,ymax=1,ytick={0,0.25,0.5,0.75,1},ylabel={fraction of runs}]
\addplot[p5blue] coordinates {(1,0) (2,0) (3,0) (5,0) (7,0) (10,0) (15,0) (20,0) (30,0) (45,0) (60,0) (90,0) (120,0) (180,0.0233) (240,0.1067) (360,0.4267) (480,0.6933) (600,0.82) (900,0.8933) (1200,0.9133) (1800,0.97) (2400,0.9967) (3600,1) (4800,1) (7200,1)};
\addplot[p5orange,dash pattern=on 4.5pt off 2.5pt] coordinates {(1,0) (2,0) (3,0) (5,0) (7,0) (10,0) (15,0) (20,0) (30,0) (45,0) (60,0) (90,0) (120,0) (180,0) (240,0) (360,0) (480,0) (600,0) (900,0.0067) (1200,0.0233) (1800,0.0467) (2400,0.05) (3600,0.0733) (4800,0.0933) (7200,0.1133)};
\nextgroupplot[title={within 0.5\%},xlabel={time (s)},xmin=1,xmax=7200,ymin=0,ymax=1,ytick={0,0.25,0.5,0.75,1}]
\addplot[p5blue] coordinates {(1,0) (2,0) (3,0) (5,0) (7,0) (10,0) (15,0) (20,0) (30,0) (45,0) (60,0) (90,0) (120,0) (180,0) (240,0) (360,0) (480,0.0133) (600,0.0533) (900,0.2) (1200,0.4133) (1800,0.7833) (2400,0.8833) (3600,0.97) (4800,0.9933) (7200,1)};
\addplot[p5orange,dash pattern=on 4.5pt off 2.5pt] coordinates {(1,0) (2,0) (3,0) (5,0) (7,0) (10,0) (15,0) (20,0) (30,0) (45,0) (60,0) (90,0) (120,0) (180,0) (240,0) (360,0) (480,0) (600,0) (900,0) (1200,0) (1800,0) (2400,0) (3600,0) (4800,0) (7200,0)};
\end{groupplot}
\coordinate (leg) at ($(group c2r2.south) + (0,-1.6cm)$);
\draw[p5blue,line width=1.2pt] ($(leg)+(-2.5,0)$) -- ($(leg)+(-1.9,0)$);
\node[anchor=west,font=\footnotesize] at ($(leg)+(-1.82,0)$) {TD-ILS};
\draw[p5orange,line width=1.2pt,dash pattern=on 4.5pt off 2.5pt] ($(leg)+(0.45,0)$) -- ($(leg)+(1.05,0)$);
\node[anchor=west,font=\footnotesize] at ($(leg)+(1.13,0)$) {TD-ACO};
\end{tikzpicture}
\caption{Anytime TD-ILS and TD-ACO convergence on Dabia2013 at $n=100$ (top) and \poryos{} at $n=1000$ (bottom), under matched budgets and seeds. The panels show mean gap and the share of runs within 1\% and 0.5\% of the BKS. References are pinned to the post-campaign store state of 2026-07-20; later improvements affect absolute Dabia2013 gaps slightly but not the strategy comparison.}
\label{fig:convergence-ils-aco}
\end{figure}
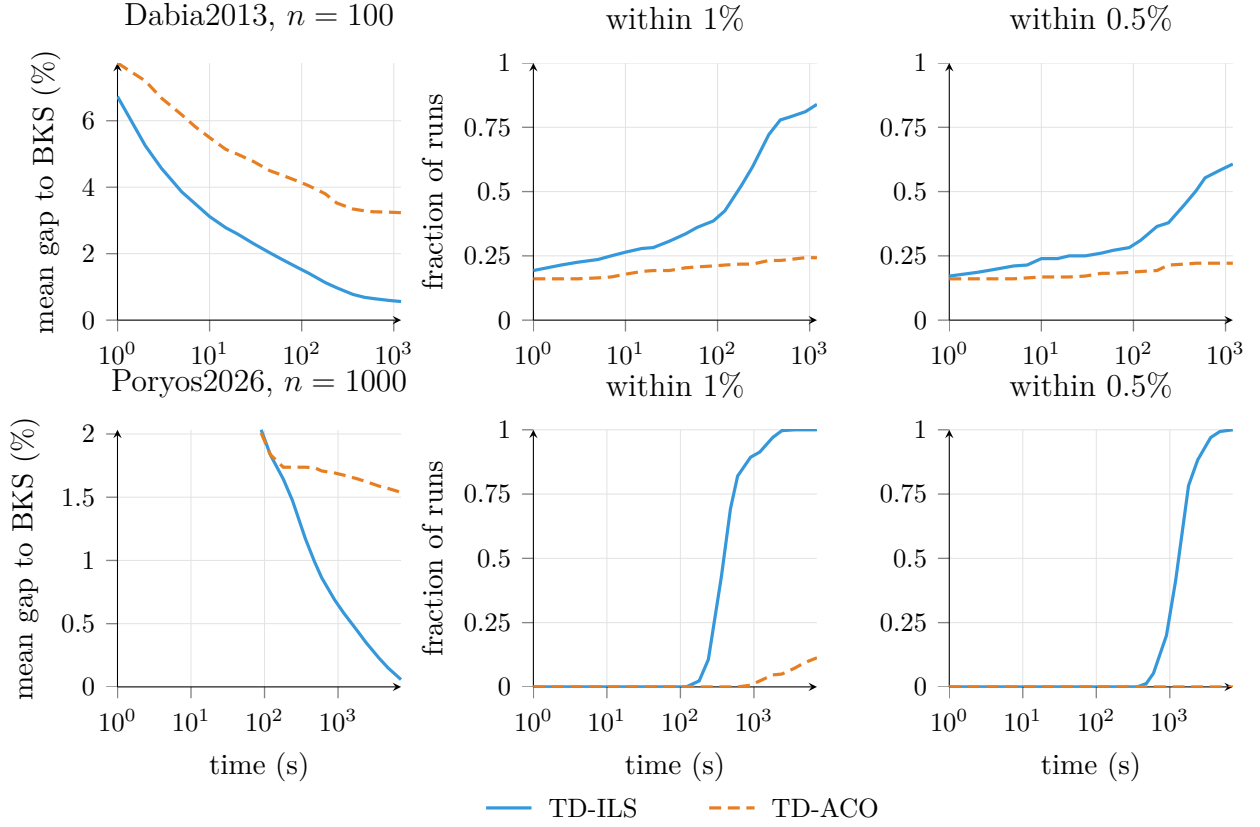

\Cref{fig:convergence-ils-aco} makes the practical difference visible on representative classic and large urban instances. ILS continues improving after ACO begins to plateau, reaches much smaller final gaps, and places far more runs within the 1\% and 0.5\% targets. The reference values are pinned to the store snapshot produced with the campaign, so later BKS improvements slightly change absolute gaps but not the paired comparison.

% ########################################################
\section{The human--AI collaboration}\label{sec:collab}
% Budget: ~2.5-3 pages. WRITTEN LAST (tone register). Methodology, not testimony.

\kayros{} had been on the lead author's roadmap for years before it existed. This section documents the working method through which it was implemented in about two weeks and states the limits of that experience report. The personal positions on AI-assisted research expressed here are those of the lead author and are endorsed by his coauthors.

\subsection{Setting and division of labor}\label{sec:collab-setting}

The engineering collaboration described in this section involved the lead author, Florian Rascoussier, then a third-year PhD student, and AI tools. The AI side is a set of frontier models, chiefly Claude Fable 5 (Anthropic), used through an agentic coding environment that can inspect and edit files, run tests, and orchestrate experiments on Grid'5000. The figures that follow cover the initial report and the first public \kayros{} release, not the later campaigns. Over roughly 14 days, the lead author estimates more than 100 hours of personal work\footnote{Active work time, excluding pauses.}, approximately $2.5 \times 10^{9}$ tokens across $10\:000$ calls, and more than \$4000 in API credits\footnote{Usage benefited substantially from Anthropic's promotional period following the release of Claude Fable 5.}. These are orders of magnitude tracked by the lead author, not audited measurements. The lead author set the research direction, decided algorithms and publication claims, governed the benchmarks and reviewed released work. The agent performed much of the implementation, debugging, experiment orchestration and initial drafting. The human coauthors retain authorship and responsibility for this report.

\subsection{Context-Oriented Programming}\label{sec:collab-cop}

The central practice was what the lead author informally calls \emph{Context-Oriented Programming} (COP): treating the curation of context, rather than the wording of an isolated prompt, as the main lever on output quality. The project draws on more than twenty experimental repositories and nearly $3\:000$ research notes. For each work stream, the relevant knowledge was distilled into design memos stating invariants and decisions, session logs recording attempts and failures, and onboarding briefs defining the mission and validation gates. Each session consumed and extended these artifacts, so project memory remained reviewable outside the model's context window. The term is informal and the method is not experimentally validated; it names the discipline that made long autonomous sessions useful in this project.

\subsection{What became feasible}\label{sec:collab-feasible}

Three concrete capabilities distinguish this project from what the lead author could have done alone, and none of them is “typing faster”. It accelerated design iteration inside inherited research code by making it affordable to instrument competing hypotheses, build reproducers and discard plausible but incorrect fixes. It made large validation campaigns routine: the agent prepared and monitored thousands of independent runs, while the lead author reviewed protocols and conclusions. Finally, it lowered the cost of replacement, allowing the engine, local search and bindings to be rewritten as plain layers instead of accumulating patches around older prototypes. The lead author had previously scoped this body of work at $2$+ years of post-doctoral effort. That comparison is a counterfactual judgment from one practitioner, not a measured productivity factor.

\subsection{Verification and accountability}\label{sec:collab-accountability}

The workflow is credible only where its evidence can be checked independently of both participants. Published routes are evaluated by an external checker, computational certificates require diversified solves and audited exact pricing, and promotions follow predefined gates (\cref{sec:certification}). These controls did not prevent invalid optimality claims from being released on their repositories. They made contradictions visible, supported diagnosis and forced retraction rather than reinterpretation. The lead author's position, argued by this report's existence, is that AI-augmented engineering does not lower the bar for verification. It raises it because verification becomes the binding constraint on what the collaboration may claim.

\subsection{Limits of this account}\label{sec:collab-limits}

This account concerns the lead author's engineering work on a single project, conducted under the scientific supervision of Romain Billot, Lina Fahed and Christine Solnon. The effort figures are self-reported, and the counterfactual (\enquote{years of solo work}) is an informed estimate, not an observation. The lead author was unusually well prepared for exactly this collaboration, with years of accumulated notes, benchmarks and prototypes for the agent to stand on, and the agent generation used here was released days before the work. Both facts limit generalization in opposite directions. Nothing here was a controlled comparison: no baseline team, no ablation of the working method, no blinded review of output quality. A proper study of AI-assisted solver engineering would need several teams, matched tasks, and independent verification budgets. Until such studies exist, accounts like this one should be read as documented experience reports, and this section has tried to be exactly that, failure included.

% ########################################################
\section{Limitations and outlook}\label{sec:limitations}

This report does not provide a head-to-head ranking against another solver. The available systems require different TD models, objectives or adapters, so a responsible comparison needs a separate comparability protocol; those results belong to the thesis and full paper. The composition engine's performance study is likewise outside the present scope. \kayros{} remains beta software, its distributed wheels target Linux x86-64, and its exact component supports the \texttt{Duration} objective but not \texttt{FleetCostDuration}. Its current certificates stop at $n=100$, with TDVRP often reaching a memory frontier before TDVRPTW.

The certificates are conditional computational results, not formal proof objects. They depend on checker-consistent route costs, standard LP and pricing tolerances, as well as labeling completeness modulo epsilon dominance. Independently verified safe dual bounds remain future work. The two recorded retractions show why this qualification matters: the repaired set has been re-derived, but continued counterexample search remains part of the protocol. Finally, benchmark names alone do not guarantee comparability. The MAMUT-routing ports have explicit conventions and canonical bytes, and values produced from another distribution must be reconciled with those conventions before they can be compared.

The two benchmark families have complementary limitations. \poryos{} uses real OSM road geometry, but demands, services, capacities, time windows and traffic are synthetic. Its controlled traffic changes costs on fixed paths rather than rerouting vehicles, and complete BKS coverage establishes feasibility baselines rather than optimality. \blauth{} uses measured-speed travel times, but its path geometries were not retained and the upstream traffic service has closed. It therefore preserves published functions without claiming full geographic regeneration. In both cases checksummed release defines the benchmark, protected from fresh extractions on evolving sources.

% ########################################################
\section{Conclusion}\label{sec:conclusion}

Time dependence changes more than an arc cost: under duration minimization, the departure time of every route becomes a decision and route evaluation becomes an optimization problem of its own. Researchers have developed methods for this setting, while users of mainstream routing tools have continued to lack a ready-made open implementation. \kayros{} fills that availability gap with checker-consistent piecewise-linear evaluation, an anytime ILS and an exact BPC component under one open installation.

The public evidence is substantial but deliberately qualified. The store carries 704 computational optimality certificates under the four-solve protocol. Earlier invalid claims were exposed by better checker-valid solutions, retracted and re-derived after repair. The resulting lesson is not that redundancy guarantees correctness, but that independent evaluation, diversified computation and continued counterexample search make computational claims easier to challenge and correct. It preserves honesty and accountability as part of the iterative research refinement process within the Open Science framework.

The benchmark contribution gives that evidence a broader setting. \blauth{} preserves measured-speed delivery instances and strong fleet-priced references up to $n=2000$; \poryos{} supplies 1,080 paired static and TD instances with retained urban geography, controlled traffic and feasibility witnesses. Together they separate two forms of realism that are often conflated: fidelity to observed travel conditions and experimental control over the factors being studied. \kayros{} improves all thirty high-effort BonnTour BKS provided with \blauth{} and ten historic Dabia2013 references, suggesting strong performance on both large realistic instances and established benchmarks. The human--AI collaboration made this breadth feasible, while human responsibility and public artifacts remain the basis on which readers can assess it.

% ########################################################
\clearpage
\appendix
\section{Per-instance literature comparisons}\label{sec:appendix-comparisons}

The main narrative reports the patterns established by the two literature comparisons. This appendix retains the per-instance values for readers who need to audit individual records; the public MAMUT-routing store remains the machine-readable authority.

\subsection{Certified improvements over Dabia2013 references}

\begin{table}[H]
\centering
% AUTO-GENERATED by scripts/gen_tables.py; certified optima read from the public BKS store,
% published values and their proven/unproven status from the canonical re-pricing of the
% Lera-Romero et al. 2020 distribution (tdvrptw-workspace/outputs/lera-bks-check-table.json).
\begin{tabular}{lrrr}
\toprule
Instance & Published best value & \kayros{} certified optimum & Improvement \\
\midrule
RC106\_50 & 11830.264 & 11756.556 & $-0.62$\% \\
R104\_100 & 18344.469 & 17949.848 & $-2.15$\% \\
R107\_100 & 18890.956 & 18704.105 & $-0.99$\% \\
R108\_100 & 17626.206 & 17400.040 & $-1.28$\% \\
R110\_100 & 18786.882 & 18607.788 & $-0.95$\% \\
RC101\_100 & 24778.294 & 24762.734 & $-0.06$\% \\
RC102\_100 & 22689.967 & 22369.064 & $-1.41$\% \\
RC103\_100 & 20680.027 & 20325.853 & $-1.71$\% \\
RC106\_100 & 21382.243 & 21233.331 & $-0.70$\% \\
RC107\_100 & 20101.592 & 19918.207 & $-0.91$\% \\
\bottomrule
\end{tabular}

\caption{The ten Dabia2013 TDVRPTW instances on which the \kayros{} computationally certified optimum lies below the best published value from \textcite{lera-romeroLinearEdgeCosts2020}. The earlier references were left unproven.}
\label{tab:improvements}
\end{table}

\subsection{Heuristic improvements over \blauth{} references}

\begin{table}[H]
\centering
% AUTO-GENERATED by scripts/gen_tables.py; current values read from the public BKS store,
% BonnTour references from the pre-Kayros snapshot reports/data/blauth2024-bonntour-reference.json.
\begin{tabular}{lrrrr}
\toprule
Instance & BonnTour (\$) & \kayros{} (\$) & Routes & Improvement \\
\midrule
berlin\_500 & 2826.49 & 2823.83 & 9 & $-0.094$\% \\
cincinnati\_500 & 2906.21 & 2896.52 & 9 & $-0.333$\% \\
kyiv\_500 & 3145.48 & 3143.59 & 10 & $-0.060$\% \\
london\_500 & 4820.82 & 4783.91 & 15 & $-0.766$\% \\
madrid\_500 & 3168.94 & 3167.92 & 10 & $-0.032$\% \\
nairobi\_500 & 3112.93 & 3072.03 & 10 & $-1.314$\% \\
new\_york\_500 & 3417.74 & 3376.62 & 11 & $-1.203$\% \\
san\_francisco\_500 & 3853.58 & 3851.57 & 12 & $-0.052$\% \\
sao\_paulo\_500 & 3909.96 & 3904.53 & 12 & $-0.139$\% \\
seattle\_500 & 3178.87 & 3167.01 & 10 & $-0.373$\% \\
berlin\_1000 & 4980.18 & 4975.37 & 16 & $-0.097$\% \\
cincinnati\_1000 & 5058.90 & 5055.28 & 16 & $-0.072$\% \\
kyiv\_1000 & 5390.71 & 5383.07 & 17 & $-0.142$\% \\
london\_1000 & 7870.40 & 7844.02 & 24 & $-0.335$\% \\
madrid\_1000 & 5420.04 & 5411.69 & 17 & $-0.154$\% \\
nairobi\_1000 & 5310.90 & 5286.61 & 17 & $-0.457$\% \\
new\_york\_1000 & 5454.59 & 5449.87 & 17 & $-0.087$\% \\
san\_francisco\_1000 & 6455.29 & 6437.71 & 20 & $-0.272$\% \\
sao\_paulo\_1000 & 6752.29 & 6734.43 & 21 & $-0.264$\% \\
seattle\_1000 & 5388.10 & 5379.64 & 17 & $-0.157$\% \\
berlin\_2000 & 8849.64 & 8844.57 & 28 & $-0.057$\% \\
cincinnati\_2000 & 8722.92 & 8717.06 & 27 & $-0.067$\% \\
kyiv\_2000 & 9567.13 & 9552.16 & 30 & $-0.157$\% \\
london\_2000 & 13273.04 & 13247.73 & 40 & $-0.191$\% \\
madrid\_2000 & 9308.40 & 9289.93 & 29 & $-0.198$\% \\
nairobi\_2000 & 8887.42 & 8867.89 & 28 & $-0.220$\% \\
new\_york\_2000 & 9560.46 & 9552.46 & 30 & $-0.084$\% \\
san\_francisco\_2000 & 11021.79 & 11010.77 & 34 & $-0.100$\% \\
sao\_paulo\_2000 & 11258.19 & 11247.72 & 35 & $-0.093$\% \\
seattle\_2000 & 9179.54 & 9173.76 & 29 & $-0.063$\% \\
\bottomrule
\end{tabular}

\caption{The thirty BonnTour references currently improved by \kayros{} (store state of 2026-08-10). Costs use the \texttt{FleetCostDuration} model of \textcite{blauthVehicleRoutingTimedependent2024}. The \emph{Routes} column gives the number of routes of the improved solution, and the final column is the relative improvement.}
\label{tab:blauth-improvements}
\end{table}

% ########################################################
\clearpage
\section*{Acknowledgements}

We thank Romain Fontaine for his help with Grid'5000, where every \kayros{} validation and certification campaign runs. We thank Gonzalo Lera-Romero for making his branch-price-and-cut solver open source \autocite{lera-romeroLinearEdgeCosts2020}: the direct inspiration for the exact component of \kayros{}. We thank Leon Lan, Niels Wouda, Wouter Kool and the other contributors to PyVRP, whose open-source framework was a major inspiration for \kayros{}' anytime heuristic layer, both for the insight that a single-trajectory iterated local search is a strong, simple and scalable anytime metaheuristic for vehicle routing and for its C++/Python architecture. We thank Thibaut Vidal, who initiated the movement of open-source vehicle-routing solvers and whose HGS-CVRP is a reference anytime metaheuristic for vehicle routing. We thank Adrien Pichon\footnote{\url{https://github.com/Anzury}; ORCID: \href{https://orcid.org/0009-0005-8630-3962}{0009-0005-8630-3962}.} for the collaboration on MAMUT-routing \autocite{pichon:hal-05629810v1}, the platform on which the benchmark catalog and best-known-solution store of \cref{sec:results} are published. We also thank Marc Sevaux and Alexandru-Liviu Olteanu, fellow members of the ANR-MAMUT project behind that catalog. \poryos{} itself was designed and generated by the lead author. We thank OpenStreetMap contributors for the geographic data on which its road networks are based \autocite{openstreetmapContributorsOpenStreetMap2026}.\footnote{Extended acknowledgements are given in the AUTHORS file of the source repository: \url{https://github.com/0nyr/kayros/blob/e30fa492466b67809cdfd1e4b945230ba0bdcb2e/AUTHORS.md}}

Experiments presented in this report were carried out using the Grid'5000 testbed, supported by a scientific interest group hosted by Inria and including CNRS, RENATER and several universities as well as other organizations\footnote{\url{https://www.grid5000.fr}}.

This work is funded by the French National Research Agency (ANR) as part of the MAMUT project, ANR-22-CE22-0016, \enquote{Machine learning And Matheuristics algorithms for Urban Transportation}.

% ########################################################
\section*{Declaration of AI use}

Generative AI tools were used throughout the preparation of this work, including frontier models such as ChatGPT 5.6 Sol (OpenAI) and Claude Fable 5, Opus 5 and Sonnet 5 (Anthropic). They contributed to the solver and benchmark-generation tooling, to experimentation and debugging, and to the drafting of this report, as documented in \cref{sec:collab}. The authors reviewed and revised all outputs from these tools and take full responsibility for their use in the final content of the work.

% ########################################################
\printbibliography

\end{document}